\documentclass{amsproc}
\usepackage{amssymb}
\usepackage{graphicx}
\usepackage{amscd}
\usepackage{amsmath}
\usepackage{showkeys}
\theoremstyle{plain}

\newtheorem{definition}{Definition}

\numberwithin{equation}{section}

\newcommand\bess{\begin{eqnarray*}}
\newcommand\eess{\end{eqnarray*}}

\newcommand\beq{\begin{equation}}
\newcommand\eeq{\end{equation}}

\def\Nset{\mbox{I\kern-.21em N}}
\def\RE{{\mbox{\rm I\kern-.21em R}}}
\def\ZZ{{\mbox{\sf Z\kern-.45em Z}}}
\def\vv{\kern.344em{\rule[.18ex]{.075em}{1.32ex}}\kern-.344em}

\def\<{\langle} \def\>{\rangle}
\def\cl{\operatorname{col}}

\begin{document}

\title[]{On inverse dynamical and spectral problems for the wave and Schr\"odinger equations on finite trees. The leaf peeling method.}
\author{S. A. Avdonin}
\address{Department of Mathematics and Statistics \\
University of Alaska Fairbanks\\
Fairbanks, AK 99775-6660, USA} \email{s.avdonin@alaska.edu}
\author{ V. S. Mikhaylov}
\address{St.Petersburg   Department   of   V.A.Steklov    Institute   of   Mathematics
of   the   Russian   Academy   of   Sciences, 7 Fontanka,
St.Petersburg  191023, Russia; and St.-Petersburg State University,
Faculty of Physics, University Embankment 7-9,
St-Petersburg 199034, Russia} \email{v.mikhaylov@spbu.ru}
\author{K. B. Nurtazina}
\address{L.N. Gumilyov Eurasian  National  University, 2, Satpayev Str., Astana, 010008, Kazakhstan}
\email{nurtazina.k@gmail.com}

\thanks{This research was supported in part by the Ministry of Education and Science of Republic of Kazakhstan, grant no. 4290/GF4. }
\thanks{Sergei Avdonin was supported by NSF grant DMS 1411564.}
\thanks{Victor Mikhaylov was supported by RFBR
14-01-00535, RFBR 14-01-31388 and NIR SPbGU 11.38.263.2014.}
\date{June 15, 2015}

\maketitle

\begin{abstract}
 Interest in inverse dynamical, spectral and scattering problems for
differential equations on graphs is motivated by possible applications to
nano-electronics and quantum waveguides and by a variety of other classical and quantum applications.
Recently a new effective leaf peeling method has been proposed by S. Avdonin and P. Kurasov \cite{AK}
for solving inverse problems on trees (graphs without cycles). It allows
recalculating efficiently the inverse data  from the
original tree to the smaller trees, `removing' leaves step by step up to the rooted
edge.
In this paper we describe the main step of the spectral
and dynamical versions of the peeling algorithm --- recalculating the inverse data for the `peeled tree'.

\end{abstract}
\section{Introduction.}

Metric graphs with defined on them differential operators or differential equations are called
quantum graphs (or differential equation networks).
There are two groups of uniqueness results concerning boundary
inverse problems for quantum trees. Brown and Weikard \cite{BW},
Yurko \cite{Y}, and  Freiling  and Yurko \cite{FY}
proved uniqueness results for trees with a priori
known topology (connectivity) and lengths of edges
using the  Titchmarsh-Weyl matrix function (TW-function)
as the inverse data.

The second group concerns inverse problems with unknown topology,
lengths of edges and potentials.  Belishev \cite{B} and
Belishev and Vakulenko \cite{BV} used the spectral data,
i.e. eigenvalues and the traces of the derivatives
of eigenfunctions at the boundary vertices
(that is equivalent to  the knowledge of the TW-function).

The paper by Avdonin and Kurasov \cite{AK} contains the most complete
results in this direction. It proves that  a
quantum tree is uniquely determined by the reduced TW-function
 associated to all except one boundary vertices.
Moreover,  the reduced response
operator (the dynamical Dirichlet-to-Neumann map) known on sufficiently large
time interval also determines a tree uniquely. The expression ``sufficiently
large''  means precisely the time interval of exact controllability of the
tree.

The most significant result of  \cite{AK}
is developing a constructive and robust procedure for the recovery tree's parameters, which became known as the {\em leaf peeling  (LP)  method}.
This method was extended to boundary inverse problems with nonstandard vertex conditions in
 \cite{AKN}, to the two-velocity wave equation in \cite{ALM,ALMR}
 and to the parabolic type  equations on trees
in  \cite{AB}.

Our procedure is  recursive --- it allows
recalculating efficiently the TW-function (and response operator)   from the
original tree to the smaller trees, `removing' leaves step by step up to the rooted
edge. Because of  its recursive  nature, this procedure contains only Volterra type equations,
and hence may be a base for developing
effective numerical algorithms.
The fact that the proposed procedure is recursive is crucial for
its numerical realization since the number of edges of graphs
arising in applications is typically very big.

The development of effective numerical algorithms for solving inverse problems
on quantum graphs is one of the goals of our research of quantum graphs. The first results in this direction for a
star graph were obtained in \cite{ABM}. The LP method combines both spectral and dynamical approaches to inverse
problems. However, the main step of the LP algorithm
--- recalculating the inverse data for the `peeled tree' --- was described in detail only in a spectral version. In this paper we describe its dynamical version, very important from theoretical
and numerical viewpoints.

Let $\Omega$ be a finite connected compact graph without cycles (a
tree). The graph consists of edges $E=\{e_1,\ldots,e_{N}\}$
connected at the vertices $V=\{v_1\ldots,v_{N+1}\}$. Every edge
$e_j\in E$ is identified with an interval $(a_{2j-1},a_{2j})$ of
the real line. The edges are connected at the vertices $v_j$ which
can be considered as equivalence classes of the edge end points
$\{a_j\}$. The boundary $\Gamma=\{\gamma_1,\ldots,\gamma_m\}$ of
$\Omega$ is a set of vertices having multiplicity one (the
exterior nodes). Since the graph under consideration is a tree,
for every $a,b\in\Omega,$ $a\not=b,$ there exists the unique path
$\pi[a,b]$ connecting these points.

\section{Spectral and dynamical problems, inverse data.}

Let $\partial w(a_j)$ denotes the derivative of $w$ at the vertex
$a_j$ taken along the corresponding edge in the direction outward
the vertex.  We associate the
following spectral problem on the graph $\Omega$ with the potential $q\in L_1(\Omega):$
\begin{eqnarray}
-\frac{d^2w}{dx^2}+qw=\lambda w,\label{Shr_eqn}\\
w\in C(\Omega),\label{Shr_cont}\\
\sum_{e_j\sim v}\partial_j w(v)=0\quad \text{for $v\in V\backslash\Gamma$}, \label{Kirh}\\
w=0\quad \text{on $\Gamma$}\label{Bound_Dir}.
\end{eqnarray}
In (\ref{Kirh}) $ \partial_j w(v)$ denotes the derivative of $w$ at the vertex $v$ taken along the edge $e_{j}$ in the direction outwards the vertex.
Also, $e_{j} \sim v$ means edge $e_{j}$ is incident to vertex $v$, and the sum is taken over all edges incident to $v$.

It is well-known fact 
that the problem (\ref{Shr_eqn})--(\ref{Bound_Dir}) has a discrete
spectrum of eigenvalues
$\lambda_1\leqslant\lambda_2 \leqslant \ldots $,
$\lambda_k\to +\infty,$ and corresponding eigenfunctions
$\phi_1,\phi_2,\ldots$ can be chosen so that
$\{\phi_k\}_{k=1}^\infty$ forms an orthonormal basis in
$\mathcal{H}:=L_{2}(\Omega)$:
\begin{eqnarray*}
(\phi_i,\phi_j)_\mathcal{H}=\int_\Omega
\phi_i(x)\phi_j(x)\,dx=\delta_{ij}\,.
\end{eqnarray*}
Set $\varkappa_k(\gamma)=\partial\phi_k(\gamma)$,
$\gamma\in\Gamma$. Let $\alpha_k$ be the $m$-dimensional column
vector defined as
$\alpha_k=\cl\left(\frac{\varkappa_k(\gamma)}{\sqrt{\lambda_k}}\right)_{\gamma\in\Gamma}.$
\begin{definition}
The set of pairs
\begin{equation}
\label{DSD} \left\{\lambda_k,\alpha_k\right\}_{k=1}^\infty
\end{equation}
is called the Dirichlet spectral data of the graph $\Omega$.
\end{definition}
The inverse spectral problem with the data given by (\ref{DSD})
was considered in \cite{B}. We will be dealing with the
Titchmarsh-Weyl function which is defined as follows. We consider
the differential equation on $\Omega$ for $\lambda\notin
\mathbb{R}:$
\begin{equation}
\label{eq} -\phi''(x)+q(x)\phi(x)=\lambda\phi(x).
\end{equation}
Let the function $\psi_i(x,\lambda)$ be the solution to (\ref{eq})
with standard conditions (\ref{Shr_cont}), (\ref{Kirh}) at
internal vertices and the following boundary conditions
\begin{equation*}
\psi_i(\gamma_i,\lambda)=1,\quad \psi_i(\gamma_j,\lambda)=0,\quad
j\not=i.
\end{equation*}
Then the entries of the Titchmarsh-Weyl matrix $\mathbf
M(\lambda)$ are defined as
\begin{equation*}
M_{ij}(\lambda)=\partial \psi_i(\gamma_j).
\end{equation*}

Let t $\psi$ be the solution to
(\ref{Shr_eqn}), (\ref{Shr_cont}) and (\ref{Kirh}) with the nonhomogeneous Dirichlet boundary conditions:
\begin{equation}
\label{SpDirF} \psi=\zeta \ \text{on} \ \Gamma,
\end{equation}
where $\zeta\in  \mathbb{R}^m.$ The Titchmarsh-Weyl $M-$matrix connects the values
of $\psi$ on the boundary and the values of its derivative
on the boundary:
\begin{equation}
\label{M_def} \partial \psi=\mathbf{M}(\lambda)\zeta \quad \text{on
$\Gamma$}.
\end{equation}

Along with the spectral we consider the dynamical problem:
\begin{eqnarray}
u_{tt}-u_{xx}+qu=0\quad\text{in $\Omega\backslash V\times [0,T]$},\label{wv_eqn}\\
u|_{t=0}=u_t|_{t=0}=0,\label{cond_1}\\
u(\cdot,t)\quad \text{satisfies $(\ref{Shr_cont})$ and $(\ref{Kirh})$ for all $t\in[0,T]$},\label{cond_2}\\
u=f\quad\text{on $\Gamma\times [0,T]$}\label{cond_3}.
\end{eqnarray}
Here $T>0$, $f=f(\gamma,t)$, $\gamma\in\Gamma$, is the Dirichlet
boundary control which belongs to
$\mathcal{F}^T_\Gamma:=L_2([0,T];\mathbb{R}^m)$. Let $u^f$ be the
solution to the problem (\ref{wv_eqn})--(\ref{cond_3}) with the
boundary control $f$. We introduce the \emph{dynamical response
operator} (the dynamical Dirichlet-to-Neumann map) by the rule
\begin{equation*}
\left(R^T\{f\}\right)(t)=\partial u^f(\cdot,t)\Bigl|_\Gamma,\quad t \in
[0,T].
\end{equation*}
The response operator has a form of  convolution:
\begin{equation*}
\left(R^T\{f\}\right)\left(t\right)=\left(\mathbf{R}*f\right)\left(t\right),
t \in [0,T].
\end{equation*}
Here the matrix-valued \emph{response function} $\mathbf{R}(t)$
is defined by the following procedure. Let the function $u_i(x,t)$
be the solution to (\ref{wv_eqn})--(\ref{cond_2}) and the
 boundary conditions
\begin{equation*}
u_i(\gamma_i,t)=\delta(t),\quad u_i(\gamma_j,t)=0,\quad j=1,\ldots
m,\,\,j\not=i.
\end{equation*}
The entries of response matrix $\mathbf{R}$ are defined by
\begin{equation*}
R_{ij}(t)=\partial u_i(\gamma_j,t).
\end{equation*}

 Connection between the spectral and dynamical data are known and were
used for studying inverse spectral and dynamical problems, see
for example \cite{AK,AMR}. Let $f\in \mathcal{F}^T_\Gamma\cap
(C_0^\infty(0,+\infty))^m$ and
$$
\widehat{f}(k):=\int_0^\infty f(t)e^{ikt}\,dt
$$
be its Fourier transform. The equations (\ref{wv_eqn}), and
(\ref{Shr_eqn}) are connected by the Fourier transformation: going
formally in (\ref{wv_eqn}) over to the Fourier transform, we
obtain (\ref{Shr_eqn}) with $\lambda=k^2$. It is not difficult to
check (see, e.g. \cite{AK, AMR}) that the response matrix-function
and Titchmarsh-Weyl matrix are connected by the same transform:
\begin{equation*}
\mathbf{M}(k^2)=\int_0^\infty \mathbf{R}(t)e^{ikt}\,dt,
\end{equation*}
where this equality is understood in a weak sense.

\section{Inverse spectral problem. Recovering the Titchmarsh-Weyl function for the peeled tree.}

Any boundary vertex of the tree can be taken as a root; therefore
without  loss of generality, we can assume that the boundary
vertex $\gamma_m$ is a root of the tree. We put $\Gamma_m=\Gamma
\setminus \{\gamma_m\}$ and consider the spectral problem
(\ref{Shr_eqn})--(\ref{Bound_Dir}) on $\Omega$. The reduced TW
matrix ${\bf M}(\lambda)=\{M_{ij}(\lambda)\}_{i,j=1}^{m-1}$
associated with boundary points from $\Gamma_m$ is constructed as
in Section 1 and serve as data for our inverse problem.

Using the procedure described in \cite{AK} we can recover the
potential on all the boundary edges (it suffices to know only the
diagonal elements of the $M-$matrix to do it). Moreover, using the
non-diagonal elements we can identify the sets of boundary
edges incident to the same internal vertex. We call these sets pre-sheaves.
More precisely, we introduce the following
\begin{definition} We consider a subgraph of $\Omega$ which is a star graph
consisting of {\sl all} edges incident to an internal vertex $v.$
This star graph is called a \emph{pre-sheaf} if it contains at
least one  boundary edge of  $\Omega.$ A pre-sheaf is called a
\emph{sheaf} if all but one its edges are the boundary edges of
$\Omega.$
\end{definition}
Following \cite{ALMR} one can extract a sheaf from all
pre-sheaves found on the previous step and proceed with the
leaf-peeling method procedure described below.

Let the  found sheaf consist of the boundary vertices
$\gamma_1,\ldots,\gamma_{m_0}$ from $\Gamma_m,$ the corresponding boundary
edges $e_1,\ldots,e_{m_0}$ and an internal edge $e_{m'_0}.$   We
identify each edge $e_i, \, i=1,\ldots,m_0$, with the interval
$[0,l_i]$ and the vertex $\gamma_{m'_0},$  the internal vertex of the
sheaf, ---  with the set of common endpoints $x=0.$ At this point
it is convenient to renumerate the edge $e_{m'_0}$ as $e_{0}$ and
the vertex  $\gamma_{m'_0}$ as  $\gamma_{0}$. Applying the techniques from
\cite{AK,ALMR} we recover the potential $q$ and lengths $l_i$ of
edges $e_i$, $i=1,\ldots,{m_0}$.

We call $\widetilde {\mathbf{M}}(\lambda)$ the reduced $M-$matrix
associated with the new graph $\widetilde \Omega=\Omega\backslash
\bigcup_{i=1}^{m_0}\{e_i,\gamma_i\}$  with boundary points
$\Gamma\backslash \bigcup_{i=1}^{m_0}\gamma_i$.

First we recalculate entries $\widetilde M_{0i}(\lambda)$,
$i=0,m_0+1,\ldots,m-1$. Let us fix $\gamma_1$, the boundary point
of the star-subgraph. Let $u$ be the solution to the problem
(\ref{eq}) with the boundary conditions
\begin{equation*}
u(\gamma_1)=1,\quad u(\gamma_j)=0,\quad j=2,\ldots,m.
\end{equation*}
We point out that on the boundary edge $e_1$ the function $u$
solves the Cauchy problem
\begin{equation}
\label{v1}
\left\{\begin{array}l
-u''(x)+q(x)u(x)=\lambda u(x), \quad x\in e_1,\\
u(l_1)=1,\,\, u'(l_1)=M_{11}(\lambda).
\end{array}
\right.
\end{equation}
On other boundary edges of the sheaf it solves the problems
\begin{equation}
\label{v3} \left\{
\begin{array}l
-u''(x)+q(x)u(x)=\lambda u(x) \\
u(l_i)=0,\,\, u'(l_i)=M_{1i}(\lambda)
\end{array}
\right.,\quad x\in e_i,\quad i=2,\ldots,m_0
\end{equation}
Since we know the potential on the edges $e_1,\ldots,e_{m_0}$, we
can solve the Cauchy problems (\ref{v1}) and (\ref{v3}), and use
the conditions (\ref{Shr_cont}), (\ref{Kirh}) at the internal
vertex $v_{0}$ to recover $u(0,\lambda)$, $u'(0,\lambda)$ -- the
values of the solution and its derivative at the ``new" boundary
edge with the ``new" boundary point $v_{0}$. Thus we obtain:
\begin{eqnarray*}
\widetilde M_{00}(\lambda)=\frac{u'(0,\lambda)}{u(0,\lambda)},\\
\widetilde
M_{0i}(\lambda)=\frac{M_{1i}(\lambda)}{u(0,\lambda)},\quad
i=m_0+1,\ldots m-1.
\end{eqnarray*}

To find $\widetilde M_{i0}(\lambda)$, $i=m_0+1,\ldots,m-1$ we fix
$\gamma_i$, $i\notin\{1,\ldots,m_0,m\}$ and consider the solution
$u$ to (\ref{eq}) with the boundary conditions
\begin{equation*}
u(\gamma_i)=1,\quad u(\gamma_j)=0,\quad j\not=i.
\end{equation*}
The function $u$ solves the Cauchy problems on the edges
$e_1,\ldots,e_{m_0}$:
\begin{equation}
\label{v5} \left\{
\begin{array}l
-u''(x)+q(x)u(x)=\lambda u(x)\\
u(\gamma_j)=0,\,\, u'(\gamma_j)=M_{ij}(\lambda)
\end{array}
\right.,\quad x\in e_j,\quad j=1,\ldots,m_0.
\end{equation}
Since we know the potential on the edges of the subgraph, we can
solve the Cauchy problems (\ref{v5}) and use the conditions at the
internal vertex $v_0$ to recover $u(0,\lambda)$, $u'(0,\lambda)$
-- the values of solution and its derivative at the ``new''
boundary edge with the ``new'' boundary point $v_0$. On the
subgraph $\widetilde\Omega$ the function $u$ solves equation
(\ref{eq}) with the boundary conditions
\begin{equation*}
u(\gamma_i)=1,\,\, u(\gamma_{0})=u(0,\lambda), \,\,
u(\gamma_j)=0,\, j=m_0+1,\ldots,
m,\,\,\gamma_j\not=\gamma_i,\,\gamma_j\not=\gamma_{0}.
\end{equation*}
Thus we obtain the equalities
\begin{eqnarray}
\label{M_ib} \widetilde
M_{i0}(\lambda)=u_x(0,\lambda)-u(0,\lambda)  \widetilde M_{00}(\lambda),
\end{eqnarray}
\begin{eqnarray}
\label{M_ibj} \widetilde
M_{ij}(\lambda)=M_{ij}(\lambda)-u(0,\lambda)  \widetilde M_{0j}(\lambda).
\end{eqnarray}
To recover all elements of the reduced matrix we need to use this
procedure for all $i,j=m_0+1,\ldots,m-1$.

Thus using the  described procedure we can identify a sheaf and
recalculate the truncated TW-matrix for the new `peeled' tree,
i.e. a tree without this sheaf. Repeating the procedure sufficient
number of times, we recover the topology of a tree and a
potential.

\section{Inverse dynamical problem. Recovering the response operator for the peeled tree.}

We assume that the boundary vertex $\gamma_m$ is a root of the tree. We
put $\Gamma_m=\Gamma \setminus \{\gamma_m\}$ and consider the dynamical
problem (\ref{wv_eqn})--(\ref{cond_3}) on $\Omega$. The reduced
response matrix function ${\bf
R}(t)=\{R_{ij}(t)\}_{i,j=1}^{m-1}$ associated with
boundary points from $\Gamma_m$ is constructed as in Section 1 and
serves as inverse data.

We use the procedure described in \cite{AK,ALM,ALMR} to recover
the potential and lengths on all boundary edges and determine all
sheaves.

Take the sheaf consisting of the boundary vertices
$\gamma_1,\ldots,\gamma_{m_0}$ from $\Gamma_m,$ the corresponding boundary
edges $e_1,\ldots,e_{m_0}$ and an internal edge $e_{m'_0}.$   We
identify each edge $e_i, \, i=1,\ldots,m_0$, with the interval
$[0,l_i]$ and the vertex $\gamma_{m'_0},$  the internal vertex of the
sheaf, ---  with the set of common endpoints $x=0.$ At this point
it is convenient to renumerate the edge $e_{m'_0}$ as $e_{0}$ and
the vertex  $\gamma_{m'_0}$ as  $\gamma_{0}$.

We call $\widetilde {\mathbf{R}}(t)$ the reduced response function
associated with the new graph $\widetilde \Omega=\Omega\backslash
\bigcup_{i=1}^{m_0} \{e_i,\gamma_i\}$ and with boundary points
$\Gamma\backslash \bigcup_{i=1}^{m_0}\gamma_i$.

To recover the entries $\widetilde R_{0i}(t)$,
$i=m_0+1,\ldots,m-1$, we consider the function $u^\delta$ to be
the solution to the problem (\ref{wv_eqn}), (\ref{cond_1}),
(\ref{cond_2}) with the boundary conditions
\begin{equation*}
u^\delta(\gamma_1,t)=\delta(t),\quad
u^\delta(\gamma_j,t)=0,\quad j=2,\ldots,m.
\end{equation*}
Since we know the potential on the edges $e_1,\ldots,e_{m_0}$, we
can evaluate $u^\delta$ on these edges by solving the wave
equations on the corresponding edges with the known boundary data.
Let us introduce the functions:
\begin{eqnarray*}
F_1(t)=\delta(t),\quad F_i(t)=0, \quad i=2,\ldots,m_0, \quad t\in \mathbb{R},\\
R_1(t)= R_{11}(t),\,\, t>0,\quad R_i(t)=\left\{ \begin{array}l 0,\quad t<l_1+l_i\\
R_{1i}(t),\quad t\geqslant l_1+l_i,
\end{array}
\right. \quad i=2\ldots,m_0.
\end{eqnarray*}
Then the function $u^\delta$ solves the following Cauchy problems
on $e_1,\ldots, e_{m_0}$:
\begin{equation*}
\left\{
\begin{array}l
u^\delta_{tt}-u^\delta_{xx}+q(x)u^\delta=0,\,\, x\in (0,l_i) \\
u^\delta(l_i,t)=F_i(t), \quad u^\delta_x(l_i,t)=R_i(t) \\
u^\delta(x,0)=0,\,\, x\in (0,l_i)
\end{array}
\right., \quad i=1\ldots,m_0
\end{equation*}
Using the compatibility conditions (\ref{cond_1}), (\ref{cond_2})
at the vertex $v_{0}$ we can find the values of $u^\delta(0,t)$
and $u^\delta_x(0,t)$ for $t>0$ at the ``new'' boundary edge $e_{0}$. We
introduce the notations
\begin{equation*}
a(t):=u^\delta(0,t),\quad A(t):=u^\delta_x(0,t),
\end{equation*}
where $a(t)=A(t)=0$ for $t<l_1$.

Let us now consider $u^f$ to be the solution to the problem
(\ref{wv_eqn}), (\ref{cond_1}), (\ref{cond_2}) with the
boundary conditions
\begin{equation*}
u^f(\gamma_1,t)=f(t),\quad u^f(\gamma_j,t)=0,\quad j=2,\ldots,m.
\end{equation*}
Due to Duhamel's principle,
$u^f(\cdot,t)=\left(u^\delta*f\right)(\cdot,t)$, and at the
``new'' boundary vertex we have:
\begin{equation*}
u^f(\gamma_{0},t)=a(t)*f(t).
\end{equation*}
By the definitions of the response matrices $\mathbf{R}$,
$\widetilde {\mathbf{R}},$ the following equalities are valid:
\begin{equation}
\label{R_m} \int_0^tR_{1i}(s)f(t-s)\,ds=\int_{0}^t\widetilde
R_{0i}(s)[a*f](t-s)\,ds, \ \ i=m_0+1,\ldots,m-1.
\end{equation}

It was proved in \cite{AB} that the response operator $R^T,$ known for sufficiently large $T,$
uniquely determines the spectral data and the Titchmarsh-Weyl matrix function. `Sufficiently large'
means precisely that $T$ is not less than the time of exact controllability:
$\,
T \geq 2 \, dist \, \{\gamma_m,\, \Gamma_m\}.\,
$ In  \cite{AB} the inverse problem for the heat equation with the Neumann-to-Dirichlet data was
studied, but the proof extends to our problem. Therefore,
we assume  below  that all sums are finite, since we may
consider the corresponding functions on finite time intervals.

 We know (see, e.g. \cite{AK}) that
$R_{1i}(s)$, $\widetilde R_{0i}(s)$, $i=m_0+1,\ldots,m-1,$ and $a$
admit the following representations (we separate regular and singular
parts):
\begin{eqnarray*}
R_{1i}(s)=r_{1i}(s)+\sum_{n\geqslant 1} \alpha_n \delta'(s-\beta_n),\quad r_{1i}|_{s\in(0,\beta_1)}=0\\
a(s)=\widetilde a(s)+\sum_{k\geqslant 1} \psi_k \delta(s-\varkappa_k),,\quad \widetilde a|_{s\in(0,\varkappa_1)}=0,\,\varkappa_1=l_1,\\
\widetilde R_{0i}(s)=\widetilde r_{0i}(s)+\sum_{l\geqslant 1}
\phi_l \delta'(s-\zeta_l),\quad \widetilde
r_{0i}|_{s\in(0,\zeta_1)}=0.
\end{eqnarray*}
In the above representations  the function $\widetilde
r_{bi}(s)$, and numbers $\psi_l$ and $\zeta_l$ are unknown and
sequences $\beta_n, \varkappa_n, \zeta_n$ are strictly increasing.
Plugging these representations to (\ref{R_m}), we obtain the following expression for the
left hand side of (\ref{R_m}):
\begin{equation}
\label{R_m1} \int_0^tR_{1i}(s)f(t-s)\,ds=\int_0^t
r_{1i}(s)f(t-s)\,ds-\sum_{n\geqslant 1}\alpha_n f'(t-\beta_n).
\end{equation}
For the right hand side of (\ref{R_m}) we have:
\begin{eqnarray}
\label{R_m2} \int_{0}^t\widetilde
R_{0i}(s)[a*f](t-s)\,ds=\int_{0}^t\widetilde
r_{0i}(s)\int_{0}^{t-s}\widetilde a(\tau)f(t-s-\tau)\,d\tau+\\
\sum_{k\geqslant 1}\psi_k\int_{0}^t\widetilde
r_{0i}(s)f(t-s-\varkappa_k)\,ds-\sum_{l\geqslant
1}\phi_l\int_{0}^{t-\zeta_l}\widetilde
a'(\tau)f(t-\zeta_l-\tau)\,d\tau-\notag\\
\sum_{l\geqslant 1}\sum_{k\geqslant 1} \phi_l\psi_k
f'(t-\zeta_l-\varkappa_k).\notag
\end{eqnarray}
Equating singular parts of the integral kernels of (\ref{R_m1})
and (\ref{R_m2}), we obtain:
\begin{equation}
\label{Sing} \sum_{n\geqslant 1}\alpha_n
\delta'(t-\beta_n)=\sum_{l\geqslant 1}\sum_{k\geqslant 1}
\phi_l\psi_k \delta'(t-\zeta_l-\varkappa_k).
\end{equation}
The equation (\ref{Sing}) allows one to recover the unknown
coefficients $\zeta_l$ and $\phi_l$. Equating the first terms in
(\ref{Sing}), we necessarily have that
\begin{equation*}
\beta_1=\zeta_1+\varkappa_1, \quad \alpha_1=\psi_1\phi_1,
\end{equation*}
and so,
\begin{equation*}
\zeta_1=\beta_1-\varkappa_1, \quad \phi_1=\frac{\alpha_1}{\psi_1},
\end{equation*}
The fact that the set $\{\alpha_1,\beta_1\},
\{\psi_1,\varkappa_1\}$ determines $\{\phi_1,\zeta_1\}$ we
represent in the following form:
\begin{equation*}
\{\alpha_k,\beta_k\}_{k=1}^{N_1},
\{\psi_k,\varkappa_k\}_{k=1}^{N_1}\Longrightarrow
\{\phi_1,\zeta_1\}, \quad m_1=1,\,\,\quad N_1=m_1.
\end{equation*}
Considering the second term in the left hand side in (\ref{Sing}),
we compare $\beta_2$ and $\zeta_1+\varkappa_2$. We get the
options:
\begin{itemize}
\item[1)] In the case $\beta_2\not=\zeta_1+\varkappa_2$ we
conclude that $\beta_2=\zeta_2+\varkappa_1$ and thus
\begin{equation*}
\zeta_2=\beta_2-\varkappa_1, \quad
\phi_2=\frac{\alpha_2}{\psi_1},\quad m_2=1.
\end{equation*}
\item[2)] In the case $\zeta_1+\varkappa_2=\beta_2$, but
$\alpha_2\not=\phi_1\psi_2$ we have that
$\zeta_2+\varkappa_1=\zeta_1+\varkappa_2=\beta_2$ and
$\alpha_2=\phi_1\psi_2+\phi_2\psi_1$, so
\begin{equation*}
\zeta_2=\beta_2-\varkappa_1, \quad
\phi_2=\frac{\alpha_2-\phi_1\psi_2}{\psi_1},\quad m_2=1.
\end{equation*}
\item[3)] In the case $\zeta_1+\varkappa_2=\beta_2$ and
$\alpha_2=\phi_1\psi_2$ we need to consider the third term in the
left hand side of (\ref{Sing}) and compare $\beta_3$ with
$\zeta_1+\varkappa_3$, using the same procedure.
\end{itemize}
Repeating this procedure sufficient number of times (say, $m_2$),
we recover $\{\phi_2,\zeta_2\}$. Suppose that to recover
$\{\phi_2,\zeta_2\}$ we used the coefficients
$\{\alpha_k,\beta_k\}_{k=1}^{N_2}$,
$\{\psi_k,\varkappa_k\}_{k=1}^{N_2};$ we express this fact in the form
\begin{equation*}
\{\alpha_k,\beta_k\}_{k=1}^{N_2},
\{\psi_k,\varkappa_k\}_{k=1}^{N_2}\Longrightarrow
\{\phi_k,\zeta_k\}_{k=1}^2, \quad N_2=N_1+m_2.
\end{equation*}

Assume that we have already recovered $p$ pairs:
\begin{equation*}
\{\alpha_k,\beta_k\}_{k=1}^{N_p},
\{\psi_k,\varkappa_k\}_{k=1}^{N_p}\Longrightarrow
\{\phi_k,\zeta_k\}_{k=1}^p.
\end{equation*}
To recover $\{\phi_{p+1},\zeta_{p+1}\}$ we need to apply the
procedure described above:  consider $N_p+1-$th term in the left hand side
of (\ref{Sing}) and compare $\beta_{N_p+1}$ with
$\zeta_p+\varkappa_{N_p+1}$ to get the options:
\begin{itemize}
\item[1)] If $\zeta_p+\varkappa_{N_p+1}\not=\beta_{N_p+1}$, then
$\beta_{N_p+1}=\zeta_{p+1}+\varkappa_{N_p}$, thus
\begin{equation*}
\zeta_{p+1}=\beta_{N_p+1}-\varkappa_{N_p},\quad
\phi_{p+1}=\frac{\alpha_{N_p+1}}{\psi_{N_p}},\quad m_{p+1}=1.
\end{equation*}
\item[2)] In the case $\beta_{N_p+1}=\zeta_p+\varkappa_{N_p+1}$,
but $\alpha_{N_p+1}\not=\phi_p\psi_{N_p+1}$, we conclude that
$\beta_{N_p+1}=\zeta_{p+1}+\varkappa_{N_p}=\zeta_p+\varkappa_{N_p+1}$
and $\alpha_{N_p+1}=\phi_p\psi_{N_p+1}+\phi_{p+1}\psi_{N_p}$. Thus
we get
\begin{equation*}
\zeta_{p+1}=\beta_{N_p+1}-\varkappa_{N_p},\quad
\phi_{p+1}=\frac{\alpha_{N_p+1}-\phi_p\psi_{N_p+1}}{\psi_{N_p}},\quad
m_{p+1}=1.
\end{equation*}
\item[3)] If $\zeta_p+\varkappa_{N_p+1}=\beta_{N_p+1}$ and
$\alpha_{N_p+1}=\phi_p\psi_{N_p+1}$ we need to compare
$\zeta_p+\varkappa_{N_p+2}$ with $\beta_{N_p+2}$ using the same
procedure.
\end{itemize}
Repeating this procedure sufficient number of times (say,
$m_{p+1}$), we recover $\{\phi_{p+1},\zeta_{p+1}\}$. We write this
in the form
\begin{equation*}
\{\alpha_k,\beta_k\}_{k=1}^{N_{p+1}},
\{\psi_k,\varkappa_k\}_{k=1}^{N_{p+1}}\Longrightarrow
\{\phi_k,\zeta_k\}_{k=1}^{p+1},\quad N_{p+1}=N_p+m_{p+1}.
\end{equation*}
The more quadruplets $\{\beta,\varkappa,\alpha,\phi\}$ we know,
the more pairs $\{\phi,\zeta\}$ we can evaluate. Certainly the
number of quadruplets at our disposal depends on the time interval
at which we know the inverse data. In the case when we know the
response function on $(0,+\infty)$, using the procedure described
above, we can recover $\{\phi_k,\zeta_k\}$ for arbitrary $k$.

Taking in (\ref{R_m1}), (\ref{R_m2}) $f(t)=\Theta(t)$ the
Heaviside function, and equating the regular parts, we arrive at
\begin{eqnarray*}
\int_0^t r_{1i}(s)\,ds=\int_{l_1}^t\widetilde
r_{0i}(s)\int_{0}^{t-s}\widetilde
a(\tau)\,d\tau\,ds\\
+\sum_{k\geqslant 1}\psi_k\int_{0}^{t-\varkappa_k}\widetilde
r_{0i}(s)\,ds-\sum_{l\geqslant
1}\phi_l\int_{0}^{t-\zeta_l}\widetilde a'(\tau)\,d\tau
\end{eqnarray*}
Everywhere below we assume that all functions are extended to the
interval $(-\infty,0)$ by zero. Differentiating the last equality
we get
\begin{equation}
\label{reg_0} r_{1i}(t)=\int_{0}^t\widetilde r_{0i}(s)\widetilde
a(t-s)\,ds+\sum_{k\geqslant 1}\psi_k\widetilde
r_{0i}(t-\varkappa_k)-\sum_{l\geqslant 1}\phi_l\widetilde
a'(t-\zeta_l),
\end{equation}
We set $s:=t-\varkappa_1$ and rewrite (\ref{reg_0}) as
\begin{eqnarray}
\label{reg_1}
r_{1i}(s+\varkappa_1)=\int_{0}^{s+\varkappa_1}\widetilde
r_{0i}(\tau)\widetilde
a(s+\varkappa_1-\tau)\,d\tau+\psi_1\widetilde
r_{0i}(s)+\\
\sum_{k\geqslant 2}\psi_k\widetilde
r_{0i}(s+\varkappa_1-\varkappa_k)-\sum_l\phi_l\widetilde
a'(s+\varkappa_1-\zeta_l).\notag
\end{eqnarray}
Let us introduce the number
\begin{equation*}
\alpha:=\min_{i\geqslant 1}(\varkappa_{i+1}-\varkappa_i).
\end{equation*}
Notice that $\alpha$ is a positive number, since we are
dealing with the problem on a finite time interval.

The integral equation (\ref{reg_1}) for the unknown function
$\widetilde r_{m_0'i}$ can be solved by steps:
\begin{itemize}
\item[1)] On the interval $(0,\zeta_1)$ we have: $\widetilde
r_{0i}(s)=0$.

\item[2)] On the interval $(\zeta_1,\zeta_1+\alpha)$ equation
(\ref{reg_1}) has the form
\begin{eqnarray}
\int_{0}^{s+\varkappa_1}\widetilde r_{0i}(\tau)\widetilde
a(s+\varkappa_1-\tau)\,d\tau+\psi_1\widetilde r_{0i}(s)=B(s),\label{reg_2}\\
B(s)=r_{1i}(s+\varkappa_1)+\sum_l\phi_l\widetilde
a'(s+\varkappa_1-\zeta_l),\notag
\end{eqnarray}
where $B(s)$ is known for $s\in (\zeta_1,\zeta_1+\alpha)$.

\item[3)] On the interval $(\zeta_1+\alpha,\zeta_1+2\alpha)$
equation (\ref{reg_1}) has form (\ref{reg_2}) where
\begin{equation*}
B(s)=-\psi_2\widetilde
r_{0i}(s+\varkappa_1-\varkappa_2)+r_{1i}(s+\varkappa_1)+\sum_{l\geqslant
1}\phi_l\widetilde a'(s+\varkappa_1-\zeta_l)
\end{equation*}
is known function on $(\zeta_1+\alpha,\zeta_1+2\alpha)$.

\item[4)] On the interval $(\zeta_1+n\alpha,\zeta_1+(n+1)\alpha)$
equation (\ref{reg_1}) has form (\ref{reg_2}) where
\begin{equation*}
B(s)=-\sum_{k=2}^{n+1}\psi_k\widetilde
r_{0i}(s+\varkappa_1-\varkappa_k)+r_{1i}(s+\varkappa_1)+\sum_{l\geqslant
1}\phi_l\widetilde a'(s+\varkappa_1-\zeta_l)
\end{equation*}
is known function on $(\zeta_1+n\alpha,\zeta_1+(n+1)\alpha)$.

\end{itemize}

To recover $\widetilde R_{00}(s)$ one need to use the following
equation:
\begin{equation*}
\int_0^t\widetilde
R_{00}(s)[a*f](t-s)\,ds=u_x^f(v_{m_0'},t)=[A*f](t).
\end{equation*}
One need to repeat the procedure described above: write down the
expansions for the $R_{00}(t)$, $A(t)$, $a(t)$ with singular and
regular parts separated, determine the singular part and afterward,
determine the regular part from corresponding integral equation.

To recover $\widetilde R_{i0}(t)$, with $i$ fixed,
$i=m_0+1,\ldots,m-1$ we consider $u^f$, the solution to the
boundary value problem (\ref{wv_eqn}) with standard conditions at
internal vertices (\ref{cond_1}), (\ref{cond_2}) and the following
boundary conditions:
\begin{equation*}
u^f(\gamma_i,t)=f(t),\quad u^f(\gamma_j,t)=0,\quad
j=1,\ldots,m,\,\, j\not=i.
\end{equation*}
Using the fact that we know the potential on the edges
$e_1,\ldots,e_{m_0}$, we can recover the solution to the problem
above by solving the Cauchy problem for the wave equations on the
corresponding edges with known boundary data. Indeed, let us introduce the
functions:
\begin{eqnarray*}
F_k(t)=0, \quad k=1,\ldots,m_0, \quad t\geqslant 0,\\
R_k(t)=[R_{ik}*f](t),\quad t\geqslant 0.
\end{eqnarray*}
The function $u^f$ solves the following Cauchy problems on
$e_1,\ldots, e_{m_0}$:
\begin{equation*}
\left\{
\begin{array}l
u^f_{tt}-u^f_{xx}+q(x)u^f=0,\quad x\in (0,l_k), \\
u^f(\gamma_i,t)=F_k(t), \quad u^f_x(\gamma_k,t)=R_k(t)
\\
u(x,0)=0,\,\, x\in (0,l_k),
\end{array}
\right.,\quad k=1,\ldots,m_0.
\end{equation*}
Using the compatibility conditions (\ref{cond_1}), (\ref{cond_2})
at the vertex $v_{0}$ we find the values of $u^f(0,t)$ and
$u^f_x(0,t)$ for $t>0$ at the new boundary vertex:
\begin{equation*}
a(t):=u^f(0,t),\quad A(t):=u^f_x(0,t).
\end{equation*}
Then using the definition of the response matrix (cf.
(\ref{M_ib})), we obtain the equalities
\begin{equation*}
\int_0^t\widetilde R_{i0}(s)f(t-s)\,ds=A(t)-[R_{00}*a](t).
\end{equation*}
These equations for $ R_{i0} $ can be analyzed using the
procedure described above: we represent $\widetilde R_{i0}(t)$, $A(t)$ and $a(t)$ as
sums of regular and singular parts, determine the singular parts of $ R_{i0} $
and then determine the regular parts from corresponding integral
equations by steps.

Concluding this section, we state that using the described  procedure
one can identify a sheaf and recalculate the truncated
response function for the new `peeled' tree, i.e. a tree without
this sheaf. Repeating the procedure sufficient number of times, we
recover the topology of the tree, potential and the lengths of the edges.

\end{document}